\documentclass[a4paper,12pt]{article}        
              
\usepackage{amsthm,amsfonts,amsmath,amssymb} 
\usepackage{mathrsfs}                        
\usepackage{dsfont}    
\usepackage{hyperref}        
\usepackage[T1]{fontenc}          
\usepackage[latin1]{inputenc}     

\usepackage[margin=7em]{geometry}    

\newcommand{\ud}{\mathrm{d}}

\newcommand{\CR}{\mathds{R}}
\newcommand{\CZ}{\mathds{Z}}

\newcommand{\depp}[2]{\frac{\partial\emph{$#1$}}{\partial{
\emph{$#2$}}}}
\newcommand{\dev}[2]{\displaystyle\frac{\ud\emph{$#1$}}{\ud\emph{$#2$}}}

\newcommand{\integral}[4]{\displaystyle\int_{\emph{$#1$}}^{\emph{$#2$}}
{\emph{$#3$}} \ \ud\emph{$#4$}}

\newcommand{\Somatorio}[2]{\displaystyle\sum\limits_{\emph{$#1$}}^
{\emph{$#2$}}}

\theoremstyle{plain} 
\newtheorem{teo}{Theorem}[section]
\newtheorem{defi}{Definition}[section]
\newtheorem{lem}{Lemma}[section]

\theoremstyle{definition}

\newtheorem{rema}{Remark}[section]

\begin{document}

\title{Darboux normal form theorem as an example of Liouville integrability theorem}
\author{Romero Solha\footnote{Departamento de Matem\'{a}tica, Universidade Federal de Minas Gerais. Postal address: Avenida Antonio Carlos, 6627 - Caixa Postal 702 -- CEP 31270-901 -- Belo Horizonte, MG. Email address: romerosolha@gmail.com } \thanks{This work was partially supported by PNPD/CAPES.}}

\date{\today}

\maketitle

%#############################################################################################

\begin{abstract}The note offers a proof of Darboux and Liouville theorems from a symplectic group action perspective. 
\end{abstract} 

%#############################################################################################

\section{Introduction}

\hspace{1.5em}Picard-Lindel\"{o}f theorem about existence and uniqueness of solutions of ordinary differential equations implies that, near a nonsingular point, a vector field can be smoothly linearised. Frobenius theorem is just a generalisation of this result for a set of vector fields that generates a Lie algebra: near a common nonsingular point the whole set of vector fields can be simultaneously smoothly linearised. The same can be said about Darboux and Liouville theorem: they are local and semilocal symplectic linearisation results. 

The original approach of presenting these classical results from a symplectic action perspective is the content of the last section of this note. From this point of view is natural to consider Darboux theorem as a consequence of Liouville theorem.

The second section presents Darboux and Liouville theorems as they usually appear in the literature. A proof of Darboux theorem (using two different techniques) and of Liouville theorem can be found in \cite{BolFom}, as well as a proof, similar to the one of this note, of theorem \ref{lemmaDL}.

For the convinience of the reader, the third section provides a compilation of definitions used in the main body of this note, its function is to fix the notation. 

Throughout this note and otherwise stated, all the objects considered will be $C^\infty$; manifolds are real, Hausdorff, paracompact, and connected; and $C^\infty(M;\CR)$ denotes the set of real-valued smooth functions over some manifold $M$.

%#############################################################################################

\section{Classical setting}

\hspace{1.5em}Darboux theorem is a normal local form result for symplectic structures stating that all symplectic manifolds look alike locally.

\begin{teo}[Darboux] Each point of a $2n$-dimensional symplectic manifold $(M,\omega)$ have a neighbourhood $V\subset M$ and coordinate functions $x_1,y_1,\dots,x_n,y_n\in C^\infty(V;\CR)$ such that $\omega\big{|}_V=\sum_{j=1}^{n}\ud x_j\wedge\ud y_j$. 
\end{teo}

The theorem asserts that a symplectic manifold $(M,\omega)$ is equivalent to the \linebreak Darboux space $(\CR^{2n},\sum_{j=1}^{n}\ud x_j\wedge\ud y_j)$, under the identification of the neighbourhood $V\subset M$ with $\CR^{2n}$. These coordinates are called Darboux coordinates and the symplectic structure $\omega$ is said to be in Darboux form when written with them, $\omega\big{|}_V=\sum_{j=1}^{n}\ud x_j\wedge\ud y_j$.

The classical Liouville theorem on the integrability of hamiltonian systems provides a semilocal normal form for the hamiltonian flow and symplectic form near a regular level set of its first integrals.

\begin{defi} An integrable system on a $2n$-dimensional symplectic manifold $(M,\omega)$ is a mapping $F=(f_1,\dots,f_n):M\to \CR^n$ such that:
\begin{itemize}
\item it is a submersion on an open dense subset of $M$;
\item its components Poisson commute amongst each other, $\{f_j,f_k\}_\omega=0$;
\item the hamiltonian vector fields generated by its components are complete\footnote{Some authors do not assume this condition, yet it holds in some cases, e.g. when the symplectic manifold is compact.}.
\end{itemize}
\end{defi}

Examples include hamiltonian systems in dimension $2$, the harmonic oscillator (in any dimension), the Kepler problem, the mathematical pendulum, the spherical pendulum, geodesic flows on surfaces of revolution, some geodesic flows on Lie groups (the free rigid body is an example of this), and various tops.  

\begin{teo}[Liouville]\label{beauty} Let $F=(f_1,\dots,f_n):M\to \CR^n$ be an integrable system on a symplectic manifold $(M,\omega)$.
\begin{itemize}
\item The hamiltonian vector fields generated by its components define an integrable (in the Sussmann \cite{SSNN} sense) distribution of the tangent bundle whose leaves are generically lagrangian, with isotropic singular leaves.
\item The connected components of the preimage of regular values (regular leaves) are homogeneous $\CR^n$ spaces; they are diffeomorphic to $\CR^{n-m}\times\mathds{T}^m$.
\item The foliation is a lagrangian fibration in a neighbourhood of each regular leaf; it defines a fibre bundle with lagrangian fibres. 
\item There are coordinates on a local trivialisation of each lagrangian leaf in which $\omega$ is in Darboux form and the flows induced by each $f_j$ are linear.  
\end{itemize}
\end{teo}

In other words, the Liouville theorem gives a description of integrable systems near the regular points of the mapping $F$.

%#############################################################################################

\section{Lie group actions on symplectic manifolds}

\hspace{1.5em}Differential geometers tend to use a vast amount of different notations for the same objects, this section is intended to clarify the ones used in this note: all the definitions and results are well cover in the literature.

\begin{defi}A smooth Lie group action on a manifold is a group homomorphism
\begin{equation}
\rho:G\to\mathrm{Diff}(M)
\end{equation}between a Lie group $G$ and $\mathrm{Diff}(M)$, the group of diffeomorphisms of $M$, such that the associated mapping from the product manifold $G\times M$ to $M$, given by 
\begin{equation}
G\times M\ni (g,p)\mapsto\rho(g)(p)\in M \ ,
\end{equation}is smooth. 
\end{defi}

The infinitesimal counterpart of a Lie group action induces a Lie algebra antihomomorphism between $(\mathfrak{g},\mathrm{ad})$, the Lie algebra of $G$, and $(\mathfrak{X}(M;\CR),[\cdot,\cdot])$, the Lie algebra\footnote{It is important to remark that if the Lie algebra structure were to be the one which coincides with the Lie algebra of the diffeomorphism group of $M$ (minus the commutator of vector fields), then a Lie group action would induce a Lie algebra homomorphism.} of smooth vector fields on $M$.

\begin{defi}The pushforward mapping of a smooth Lie group action \linebreak $\rho:G\to\mathrm{Diff}(M)$ at the identity $e\in G$ is denoted by $\rho_{*_e}:\mathfrak{g}\to\mathfrak{X}(M;\CR)$ and defined by
\begin{equation}
(\rho_{*_e}(x+cy)(f))(p):=\dev{}{t}f\circ\rho\circ\exp(tx)(p)\Big{|}_{t=0}+c\left(\dev{}{t}f\circ\rho\circ\exp(ty)(p)\Big{|}_{t=0}\right) \ ,
\end{equation}where $p\in M$, $f\in C^\infty(M;\CR)$, $x,y\in\mathfrak{g}$ and $c\in\CR$. Wherefore, the flow at time $t\in\CR$ of $\rho_{*_e}(x)\in\mathfrak{X}(M;\CR)$ is $\rho\circ\exp(tx)\in\mathrm{Diff}(M)$, and $\rho_{*_e}(x)$ is a complete vector field. 
\end{defi}

\begin{rema}
The $C^\infty(M;\CR)$-module of smooth vector fields of a manifold $M$ will be denoted by $\mathfrak{X}(M;\CR)$ when smooth vector fields are seen as derivations of the \linebreak commutative algebra $C^\infty(M;\CR)$. The $C^\infty(M;\CR)$-module $\Omega^1(M;\CR)$ of differential \linebreak one forms on $M$ is by definition the dual $C^\infty(M;\CR)$-module of $\mathfrak{X}(M;\CR)$, whilst 
$\Omega^k(M;\CR):=\mathrm{Hom}_{C^\infty(M;\CR)}(\wedge^k_{C^\infty(M;\CR)}\mathfrak{X}(M;\CR);C^\infty(M;\CR))$. 
\end{rema}

The nondegeneracy of a symplectic form $\omega\in\Omega^2(M;\CR)$ of a symplectic manifold $(M,\omega)$ induces two particular Lie subalgebras of $(\mathfrak{X}(M;\CR),[\cdot,\cdot])$, the vector space of symplectic and hamiltonian vector fields, respectively denoted by $\mathfrak{X}^{\omega}(M;\CR)$ and $\mathfrak{X}^{ham}_{\omega}(M;\CR)$. They are, respectively, isomorphic to the space of closed and exact $1$-forms on $M$.

\begin{defi}A vector field $X\in\mathfrak{X}(M;\CR)$ of a symplectic manifold $(M,\omega)$ is a symplectic vector field if $\imath_X\omega$ is closed. In the particular case where $\imath_X\omega$ is exact, $X$ is said to be a hamiltonian vector field; and a function satisfying $\imath_{X}\omega=-\ud f$ is called a hamiltonian function for $X$.
\end{defi}

The mapping ${grad}_\omega:C^\infty(M;\CR)\to\mathfrak{X}(M;\CR)$ associates to each function \linebreak $f\in C^\infty(M;\CR)$ a hamiltonian vector field, ${grad}_\omega(f)=X\in\mathfrak{X}(M;\CR)$, via the equation 
\begin{equation}
\imath_{X}\omega=-\ud f \ ,
\end{equation}which has a unique solution due to the nondegeneracy of the symplectic form.

There is not only special vector fields on a symplectic manifold $(M,\omega)$, there are also distinct types of submanifolds. If the symplectic form $\omega$ vanishes when restricted to vector fields tangent to a submanifold, this submanifold is called isotropic, and it is lagrangian when its dimension is half of the dimension of $M$.

A symplectic structure endows the space of smooth functions with a Lie algebra structure satisfying a Leibniz rule: a Poisson structure.

\begin{defi}The Poisson bracket of two functions $f_1,f_2\in C^\infty(M;\CR)$ on a symplectic manifold $(M,\omega)$ is the function defined by 
\begin{equation}
\{f_1,f_2\}_\omega:=\omega({grad}_\omega(f_1),{grad}_\omega(f_2)) \ .
\end{equation} 
\end{defi}

When a Lie group acts on a symplectic manifold preserving the symplectic structure one says that it acts symplectically. It might also happen that the infinitesimal action of a Lie group acts on a symplectic manifold via hamiltonian vector fields.

\begin{defi}A smooth Lie group action $\rho:G\to\mathrm{Diff}(M)$ on a symplectic \linebreak manifold $(M,\omega)$ is said to be symplectic when $\rho(g)^*(\omega)=\omega$ for all $g\in G$. So the image of the action, $\rho(G)\subset\mathrm{Diff}(M)$, is a subgroup of the group of symplectic diffeomorphisms of $(M,\omega)$. It is said to be hamiltonian when $\rho_{*_e}(\mathfrak{g})\subset\mathfrak{X}^{ham}_{\omega}(M;\CR)$. 
\end{defi}

A special instance of a hamiltonian action of a Lie group occurs when the infinitesimal action respects the Lie algebra structure of smooth functions provided by the Poisson bracket, some authors call these actions Poisson actions.   

\begin{defi}Let $\rho:G\to\mathrm{Diff}(M)$ be a hamiltonian action, a comomentum \linebreak mapping for this action is a $\CR$-linear mapping $\mu^*:\mathfrak{g}\to C^\infty(M;\CR)$ satisfying \linebreak $grad_\omega\circ\mu^*=\rho_{*_e}$. If in addition $\mu^*$ is a Lie algebra antihomomorphism between $(\mathfrak{g},\mathrm{ad})$ and $(C^\infty(M;\CR),\{\cdot,\cdot\}_\omega)$, then it is called an equivariant comomentum mapping. 
\end{defi}

There is also an equivalent notion dual to the one of a comomentum mapping which is usually more discussed in the literature.

\begin{defi}A momentum mapping for a hamiltonian action $\rho:G\to\mathrm{Diff}(M)$ is a mapping $\mu:M\to\mathfrak{g}^*$ such that, for each $x\in\mathfrak{g}$ and $p\in M$, the function defined by $f(x)(p):=\mu(p)(x)$ is a hamiltonian function for $\rho_{*_e}(x)$. And if $\mu\circ\rho(g)^{-1}=Ad_g^*\circ\mu$ for all $g\in G$, then it is called an equivariant momentum mapping.  
\end{defi}

It is easy to check that the existence of a momentum mapping is equivalent to the existence of a comomentum mapping, and the equivariance of $\mu$ is equivalent to $\mu^*$ be a Lie algebra antihomomorphism.

%#############################################################################################

\section{Hamiltonian action approach}

\hspace{1.5em}Integrable systems form a particular class of examples of hamiltonian $\CR^n$-actions admitting equivariant comomentum mappings.

\begin{defi} An integrable system \`{a} la Liouville on a symplectic manifold $(M,\omega)$ is a hamiltonian $\CR^n$-action $\rho:\CR^n\to\mathrm{Diff}(M)$, whose stabiliser subgroups are discrete over an open dense subset of $M$, together with an equivariant comomentum mapping $\mu^*:\CR^n\to C^\infty(M;\CR)$.
\end{defi}

Supposing that $\rho:G\to\mathrm{Diff}(M)$ is an action of the additive Lie group $G=\CR^n$, for each basis of its Lie algebra $v_1,\dots,v_n\in\mathfrak{g}=\CR^n$ one can associate an integrable distribution $\mathcal{P}:=\langle \rho_{*_e}(v_1),\dots,\rho_{*_e}(v_n) \rangle_{C^\infty(M;\CR)}\subset\mathfrak{X}(M;\CR)$ ---the vector fields are complete, and they all commute amongst each other because $\rho_{*_e}$ is a Lie antihomomorphism. The orbits (or integral leaves) of this distribution passing through a point $p\in M$ are diffeomorphic to the quotient of $\CR^n$ by the stabiliser subgroup $G_p:=\{g\in\CR^n \ ; \ \rho(g)(p)=p \}$. The action is actually given by the joint flow of the vector fields $\rho_{*_e}(v_j)$. In an open set $V\subset M$ where each $p,q\in V$ satisfies $G_p\cong G_q$ one has a fibre bundle. 

In case this action is hamiltonian, each $\rho_{*_e}(v_j)$ belongs to $\mathfrak{X}^{ham}_{\omega}(M;\CR)$, and a comomentum mapping $\mu^*:\mathfrak{g}\to C^\infty(M;\CR)$ can be linearly defined by $\mu^*(v_j):=f_j$, with $f_j\in C^\infty(M;\CR)$ an arbitrary hamiltonian function for $\rho_{*_e}(v_j)$. The Lie algebra is abelian and $\mu^*$ is a Lie antihomomorphism if and only if $\{f_j,f_k\}_\omega=0$ for all $j,k$. Thus, in order to have an equivariant comomentum mapping from this construction, the choice of hamiltonian functions must be such that $\{f_j,f_k\}_\omega=0$ for all $j,k$ (their Poisson bracket is always a constant\footnote{This constant is actually a $2$-cocycle in the Lie algebra cohomology of $\mathfrak{g}$ with values in $\CR$. For abelian Lie algebras $\mathfrak{g}$, any $2$-cochain is a $2$-cocycle and $2$-cochains are simply skewsymmetric bilinear mappings from $\mathfrak{g}\oplus\mathfrak{g}$ to $\CR$.}, but not necessarily zero).  

It is clear now that the momentum mapping associated to an equivariant $\mu^*$, if denoted by $F:M\to\CR^n$ (after the identification $\mathfrak{g}^*\cong\CR^n$), is an integrable system when the stabiliser subgroups are discrete over an open dense subset of $M$. The hamiltonian vectors $\rho_{*_e}(v_1),\dots,\rho_{*_e}(v_n)$ provide a basis for the tangent space of an orbit at any of its points, and $\omega(\rho_{*_e}(v_j),\rho_{*_e}(v_k))=\{f_j,f_k\}_\omega=0$; therefore, each orbit passing through $p\in M$ is an isotropic submanifold given by the connected components of the preimage by the momentum mapping of $F(p)\in\CR^n$. 

Under this hypothesis regular orbits (the ones associated with discrete stabiliser subgroups which are the connected components of the preimage of regular values of the momentum mapping) are diffeomorphic to $\CR^{n-m}\times\mathds{T}^m$, where the stabiliser subgroups are isomorphic to $a_p\cdot\CZ^m$, with $a_p\in\CR^m$ being the periods of the periodic hamiltonian vector fields $\rho_{*_e}(v_{j_1}),\dots,\rho_{*_e}(v_{j_m})$ passing through $p\in M$. 
       
In conclusion, integrable systems induce a foliation on the symplectic manifold whose leaves are generically lagrangian (with isotropic singular leaves) and diffeomorphic to $\CR^{n-m}\times\mathds{T}^m$, and near each regular leaf the foliation is a lagrangian fibration. 

Thus, the missing piece from Liouville theorem is the symplectic linearisation of the hamiltonian action near a regular orbit. One needs to prove a technical lemma ---Poincar\'{e} lemma for regular foliations--- before proving the linearisation of the hamiltonian action.

\begin{lem}\label{lemmaPL} Let $\alpha\in\Omega^k(M;\CR)$ be a given closed $k$-form, $\ud\alpha=0$, whose restriction to an integrable distribution of constant rank $\mathcal{P}\subset\mathfrak{X}(M;\CR)$ vanishes. Then, for each trivialising neighbourhood $A\subset M$ of the regular foliation defined by $\mathcal{P}$, there exists a $\beta\in\Omega^{k-1}(A;\CR)$ such that $\beta$ vanishes when restricted to $\mathcal{P}$ and $\alpha=\ud\beta$ on $A$. 
\end{lem}
\textit{Proof:} Frobenius theorem (or Sussmann's theorem \cite{SSNN} for constant rank distributions) implies that near each point of $M$ there exists a neighbourhood $A\subset M$, the so-called trivialising neighbourhood, diffeomorphic to $N\times\CR^{rank(\mathcal{P})}$ with $N$ being the orbit (or leaf) of the distribution $\mathcal{P}$ passing through the point. 

Let $X\in\mathfrak{X}(A;\CR)$ be the vector field whose flow at time $t\in\CR$ is the diffeomorphism $\exp(tX)\in\mathrm{Diff}(A)$ defined by $\exp(tX)(p):=(q,v-tv)$, where, under the identification $A\cong N\times\CR^{rank(\mathcal{P})}$, $p=(q,v)\in A$ with $q\in N$ and $v\in\CR^{rank(\mathcal{P})}$. 

The homotopy formula for $X\in\mathfrak{X}(A;\CR)$ applied to $\alpha\in\Omega^k(M;\CR)$ gives, on the trivialising neighbourhood $A$, 
\begin{equation}
\alpha=\ud(-H_X(\alpha))-\exp(X)^*(\alpha) \ ,
\end{equation}with
\begin{equation}
H_X(\alpha)=\integral{0}{1}{\exp(tX)^*(\imath_X\alpha)}{t} \ .
\end{equation}The reader will notice that both $N\times\{0\}\subset A$ and $\mathcal{P}$ are invariant by the flow of $X$, and that at time one any point of $A$ is mapped into $N\times\{0\}$. These properties imply that $\exp(X)^*(\alpha)=0$ and that $\beta=-H_X(\alpha)$ is a $(k-1)$-form on $A$ vanishing when restricted $\mathcal{P}$ satisfying $\alpha=\ud\beta$.\hfill $\blacksquare$\par\vspace{0.5em}

\begin{teo}\label{teoLS} The hamiltonian $\CR^n$-action of an integrable system on $(M,\omega)$ can be symplectically linearised near each of its regular orbits. 
\end{teo}
\textit{Proof:} Near each regular orbit there exists a trivial fibre bundle structure with the momentum mapping as the projection. In this local trivialisation of this lagrangian fibration, the coordinates of the basis are given by the functions $f_1,\dots,f_n\in C^\infty(M;\CR)$ and the fibres are covered by coordinate functions $y_1,\dots,y_n\in C^\infty(A;\CR)$, with $A\cong (\CR^{n-m}\times\mathds{T}^m)\times\CR^n$ and the $m$ functions $y_{n-m+1},\dots,y_n$ periodic with periods given by $a_p\in\CR^m$.

Applying lemma \ref{lemmaPL} to the symplectic form $\omega$ on the trivial fibre bundle near a regular orbit, with the integrable distribution $\mathcal{P}:=\langle \rho_{*_e}(v_1),\dots,\rho_{*_e}(v_n) \rangle_{C^\infty(M;\CR)}$, one has $\theta\in\Omega^1(A;\CR)$ satisfying $\omega\big{|}_A=\ud\theta$ and $\theta\big{|}_{\mathcal{P}}=0$. 

Since $\theta\big{|}_{\mathcal{P}}=0$, for each $X_j:=\rho_{*_e}(v_j)$ it holds $\imath_{X_j}\theta=0$, and because $\omega\big{|}_A=\ud\theta$ one has $\imath_{X_j}(\ud\theta)=-\ud f_j$; Thus, the Lie derivative of $\theta$ with respect to $X_j$ is
\begin{equation}
\pounds_{X_j}(\theta)=\imath_{X_j}(\ud\theta)+\ud (\imath_{X_j}\theta)=-\ud f_j \ . 
\end{equation}
  
The condition $\theta\big{|}_{\mathcal{P}}=0$ also implies that $\theta=-\sum_{k=1}^{n}\theta_k\ud f_k$ and the previous equation reads
\begin{eqnarray}
-\ud f_j=\pounds_{X_j}(\theta)&=&\pounds_{X_j}\left(-\Somatorio{k=1}{n}\theta_k\ud f_k\right)=-\Somatorio{k=1}{n}\pounds_{X_j}(\theta_k\ud f_k) \nonumber \\
&=&-\Somatorio{k=1}{n}\left(X_j(\theta_k)\ud f_k+\theta_k\pounds_{X_j}(\ud f_k)\right) \nonumber \\
&=&-\Somatorio{k=1}{n}\left(X_j(\theta_k)\ud f_k+\theta_k\imath_{X_j}(\ud\circ\ud f_k)+\theta_k\ud (\imath_{X_j}\ud f_k)\right) \nonumber \\
&=&-\Somatorio{k=1}{n}X_j(\theta_k)\ud f_k \ , 
\end{eqnarray}yielding $X_j(\theta_k)=\delta_{jk}$. 

The nondegeneracy of $\omega$ actually implies that $\ud\theta_1,\dots,\ud\theta_n$ is a basis: 
\begin{equation}
\omega\big{|}_A=\Somatorio{k=1}{n}\ud f_k\wedge\ud\theta_k \ .
\end{equation}Thus, the mapping defined by $(f_1,\dots,f_n,y_1,\dots,y_n)\mapsto (f_1,\dots,f_n,\theta_1,\dots,\theta_n)$ is a diffeomorphism of $A$.

The theorem is proved by now: in the coordinates $(f_1,\dots,f_n,\theta_1,\dots,\theta_n)$ the symplectic form is just the Darboux form on $A$ and the hamiltonian action is linear, i.e. it is given by 
\begin{equation}
(f_1,\dots,f_n,\theta_1,\dots,\theta_n)\mapsto (f_1,\dots,f_n,\theta_1+t_1,\dots,\theta_n+t_n) \ ,
\end{equation}where $(t_1,\dots,t_n)\in\CR^n$, because $X_j(\theta_k)=\delta_{jk}$.\hfill $\blacksquare$\par\vspace{0.5em} 

This theorem also holds true near each nondegenerate compact leaf \cite{eli1,eli2,Mi}. This generalisation is nontrivial by the simple reason that (in general) there is no Poincar\'{e} lemma for singular foliations, even in this particular case of a foliation coming from an integrable system with nondegenerate type of singularities \cite{MiSolha}.

The next theorem is an existence theorem for integrable systems near any point of a symplectic manifold.

\begin{teo}\label{lemmaDL}Near any of its points, a symplectic manifold $(M,\omega)$ of dimension $2n$ always admits a free hamiltonian $\CR^n$-action, together with an equivariant comomentum mapping.
\end{teo}
\textit{Proof:} Let $f_1\in C^\infty(M;\CR)$ be a function whose hamiltonian vector field, denoted by $X_1\in\mathfrak{X}(M;\CR)$, does not vanish at the point $p\in M$. One can always construct such a function: indeed, in a neighbourhood $V_0$ of the point, $p\in V_0\subset M$, one can use coordinates $z_1,\dots,z_{2n}\in C^\infty(V_0;\CR)$ such that $z_j(p)=0$ to define the function $(z_1+1)^2\in C^\infty(V_0;\CR)$, which extends trivially to all of $M$ by the use of bump functions, and the nondegeneracy of $\omega$ guarantees that its hamiltonian vector field has the desired property.

Thus, the vector field $X_1\in\mathfrak{X}(M;\CR)$ can be linearised near $p\in M$, i.e. there exist a neighbourhood $V_1$ containing the point, $p\in V_1\subset M$, and coordinates $x_1,\dots,x_{2n-1},y_1\in C^\infty(V_1;\CR)$ defined on it, where $X_1$can be written as $\depp{}{y_1}$. Since $X_1(f_1)=\{f_1,f_1\}_\omega=0$, this implies that $f_1$ is independent of the coordinate $y_1$. 

Also, there always exists another function $f_2\in C^\infty(M;\CR)$ satisfying, near $p\in M$, both $\{f_1,f_2\}_\omega=0$ and that its hamiltonian vector field $X_2\in\mathfrak{X}(M;\CR)$ is linearly independent of $X_1$, as long as $n$ is bigger than $1$ ---it is not difficult to see, using the local coordinates, that these conditions define a homogeneous underdetermined system of linear equations.

Now, by Frobenius theorem, the hamiltonian vector fields $X_1$ and $X_2$ define a (regular) foliation on $V_1$; therefore, there exists a possibly smaller neighbourhood, $V_2$, of $p\in M$, with coordinates $x_1,\dots,x_{2n-2},y_1,y_2\in C^\infty(V_2;\CR)$ defined on it, where $X_1$ and $X_2$ can be written as $\depp{}{y_1}$ and $\depp{}{y_2}$, respectively, and both functions $f_1$ and $f_2$ are independent of $y_1$ and $y_2$. Repeating the argument of the previous paragraph, one can find a third function, $f_3\in C^\infty(M;\CR)$, such that $X_1$, $X_2$, and $X_3$ are linearly independent near $p\in M$ and $\{f_1,f_3\}_\omega=\{f_2,f_3\}_\omega=0$, with $X_3\in\mathfrak{X}(M;\CR)$ its hamiltonian vector field. 

This reasoning works as long as the number of functions $f_j$'s is not bigger than $n$, otherwise one would reach a homogeneous system of linear equations whose solution is only the trivial one. 

Thus, for a given point of $(M,\omega)$ there exist a neighbourhood and $n$ functions defined on it providing a free hamiltonian $\CR^n$-action, together with an equivariant comomentum mapping. The action is given by the joint flow of the hamiltonian vector fields $X_1,\dots,X_n\in\mathfrak{X}(M;\CR)$ restricted to the common neighbourhood where they do not vanish, and the comomentum mapping is the linear mapping that maps a fixed basis of $\CR^n$ onto the set of functions $f_1,\dots,f_n\in C^\infty(M;\CR)$.\hfill $\blacksquare$\par\vspace{0.5em}

The reader will recognise that theorem \ref{teoLS} applied to the free hamiltonian \linebreak $\CR^n$-action of theorem \ref{lemmaDL} is just Darboux theorem in disguise, as a free hamiltonian $\CR^n$-action together with an equivariant comomentum mapping is just an integrable system \`{a} la Liouville.

%############################################################################################################################

\subsection{Equivariant normal forms and noncommutativity}

\hspace{1.5em}Let $\rho:G\to\mathrm{Diff}(M)$ be a symplectic action of a Lie group $G$ on a symplectic manifold $(M,\omega)$. An equivariant diffeomorphism with respect to $\rho$ is an element of the centre of $\rho(G)\subset\mathrm{Diff}(M)$.

\begin{teo} Let $\omega_0,\omega_1\in\Omega^2(M;\CR)$ be two symplectic forms on $M$ such that a compact Lie group $G$ acts on in a symplectic fashion, via $\rho:G\to\mathrm{Diff}(M)$, for both symplectic structures. If $N\subset M$ is a submanifold invariant by the group action, $\rho(g)(N)\subset N$ for all $g\in G$, where both symplectic structures coincide, $\omega_0\big{|}_N=\omega_1\big{|}_N$, then there exist an invariant neighbourhood $V$ of $N$, $N\subset V\subset M$ and $\rho(g)(V)\subset V$ for all $g\in G$, and an equivariant diffeomorphism $\phi\in\mathrm{Diff}(M)$ satisfying $\phi\big{|}_N=\mathrm{Id}_N$ and $\phi^*(\omega_1)=\omega_0$.  
\end{teo}

This is the equivariant version of Darboux and Weinstein's theorems when there is a symplectic action, and the proof of the linearisation theorem (theorem \ref{teoLS}) can be adapted to this situation. It would not be difficult, then, to understand the lack of a unique normal form \cite{Mel}. For the abelian case the Poincar\'{e} lemma (lemma \ref{lemmaPL}), indeed, provides a primitive for the symplectic form, whilst in the nonabelian situation the lemma can only guarantee a primitive for the difference between the symplectic form $\omega$ near the invariant submanifold $N\subset M$ and the constant symplectic form that coincides with $\omega$ at $N$. 

This approach can probably provide another proof for Mishchenko and Fomenko's conjecture on noncommutative integrable systems in the smooth category \cite{Fome,BolJa}. Frobenius theorem gives a description near regular orbits and the linearisation result, which follows from the Poincar\'{e} lemma (lemma \ref{lemmaPL}) ---and the main ingredient for its proof is Frobenius theorem. One only needs to apply these results to the hamiltonian distribution constructed from a Lie subalgebra of noncommuting first integrals.

\end{document}